\documentclass[11pt]{amsart}
\usepackage[T1]{fontenc}
\usepackage{amsmath,amssymb,amsthm, fullpage, amsaddr}
\usepackage{tikz}
\usepackage[colorlinks=true,linkcolor=blue,citecolor=blue,urlcolor=blue]{hyperref}

\newtheorem{theorem}{Theorem}[section]
\newtheorem{lemma}[theorem]{Lemma}

\newtheorem{conjecture}{Conjecture}
\theoremstyle{remark}
\newtheorem{remark}[theorem]{Remark}

\newcommand{\cs}{\chi'_{s}}
\newcommand{\csf}{\chi'_{s,f}}
\newcommand{\cD}{\mathcal{D}}
\newcommand{\cO}{\mathcal{O}}

\title[Fractional and asymptotic Brualdi--Quinn]{The Brualdi--Quinn  conjecture on strong edge-coloring  holds fractionally and asymptotically}

\author{Rong Luo$^1$, Gexin Yu$^2$, Xiangqian Zhou$^3$}
\address{$^1$Department of Mathematics, West Virginia University, Morgantown, WV 26506\\
$^2$Department of Mathematics, William \& Mary, Williamsburg, VA 23185\\
$^3$Department of Mathematics and Statistics, Wright State University, Dayton, Ohio 45435}

\email {rluo@mail.wvu.edu, gyu@wm.edu, xiangqian.zhou@wright.edu}

\begin{document}

\begin{abstract}
Let $G$ be a bipartite graph with parts $A$ and $B$, and let $\Delta_A$ and $\Delta_B$ be the maximum degrees of the vertices in $A$ and in $B$, respectively. Brualdi and Quinn conjectured that the strong chromatic index of $G$ is at most $\Delta_A\Delta_B$. We prove that the fractional strong chromatic index of $G$ is at most $\Delta_A\Delta_B$, and that the conjecture holds asymptotically: $\cs(G)\le(1+o(1))\Delta_A\Delta_B$ as $\max\{\Delta_A,\Delta_B\}\to\infty$. The error term is $12(\ln\Delta)^{1/2}\Delta^{-1/4}$ with $\Delta=\max\{\Delta_A,\Delta_B\}$.  Both results rest on a smart partition of edges and Lov\'asz Local Lemma.


\end{abstract}

\maketitle

\section{Introduction}

All graphs in this paper are finite and simple. A \emph{strong edge-coloring} of a graph $G$ is an assignment of colors to the edges of $G$ in which every color class is an induced matching; equivalently, two edges receive distinct colors whenever they share an endpoint or are joined by an edge. The \emph{strong chromatic index} $\cs(G)$ is the minimum number of colors in a strong edge-coloring of $G$. The notion was introduced by Fouquet and Jolivet~\cite{FJ}. For a graph of maximum degree $\Delta$, a greedy argument gives $\cs(G)\le2\Delta^2-2\Delta+1$. Erd\H{o}s and Ne\v{s}et\v{r}il~\cite{EN} conjectured that
\[
\cs(G)\le
\begin{cases}
\frac54\Delta^2, & \text{if $\Delta$ is even},\\[2mm]
\frac14(5\Delta^2-2\Delta+1), & \text{if $\Delta$ is odd}.
\end{cases}
\]

Molloy and Reed~\cite{MR} established the first asymptotic improvement over the greedy bound on $\cs(G)$ where they proved that $\cs(G)\le1.998\Delta^2$ for sufficiently large $\Delta$. This was successively improved to
$1.93\Delta^2$ by Bruhn and Joos~\cite{BJ}, to $1.835\Delta^2$ by
Bonamy, Perrett and Postle~\cite{BPP}, and to $1.772\Delta^2$ by
Hurley, de Joannis de Verclos and Kang~\cite{HJK}. More recently,
Davey, Hurley, de Joannis de Verclos, Kang and Volec~\cite{DHJKV}
obtained the bound $1.73\Delta^2$ by way of flag algebra.

For bipartite graphs, Faudree, Gy\'arf\'as, Schelp and Tuza~\cite{FGST} proposed the following conjecture.

\begin{conjecture}[Faudree, Gy\'arf\'as, Schelp and Tuza~\cite{FGST}]\label{conj:FGST}
Every bipartite graph $G$ with maximum degree $\Delta(G)$ satisfies $\cs(G)\le\Delta(G)^2$.
\end{conjecture}

Brualdi and Quinn~\cite{BQM} proposed a stronger, asymmetric version. Let
$G$ be bipartite with parts $A$ and $B$, and let 
\[
\Delta_A=\max_{x\in A}d(x),\qquad \Delta_B=\max_{b\in B}d(b).
\]

\begin{conjecture}[Brualdi and Quinn~\cite{BQM}]\label{conj:BQM}
Every bipartite graph $G$ satisfies $\cs(G)\le\Delta_A\Delta_B$.
\end{conjecture}

Both conjectures are tight: if $|A|=\Delta_B$ and $|B|=\Delta_A$, the complete bipartite graph on $A\cup B$ has $\Delta_A\Delta_B$ edges and all edges need distinct colors.

Steger and Yu~\cite{SY} proved Conjecture~\ref{conj:FGST} for $\Delta(G)=3$. Conjecture~\ref{conj:BQM} is confirmed when one of $\Delta_A,\Delta_B$ is at most $2$ (Nakprasit~\cite{N}) or at most $3$ (Huang, Yu and Zhou~\cite{HYZ}). For large degrees, Hao, Yang and Yu~\cite{HYY} proved $\cs(G)\le1.676\,\Delta_A\Delta_B$ when $\Delta_A \Delta_B$ is sufficiently large. Davey, Hurley, de Joannis de Verclos, Kang and Volec~\cite{DHJKV} proved $\cs(G)\le1.6255\,\Delta(G)^2$ for bipartite $G$ with $\Delta(G)$ sufficiently large, $\cs(G)\le1.6633\,\Delta_A\Delta_B$ when $\Delta_B/\Delta_A$ is a fixed rational number and $\Delta_A$ is sufficiently large, and Conjecture~2 holds asymptotically almost surely for random bipartite graphs whose two parts have comparable sizes. Before~\cite{DHJKV, HYY}, no general improvement over the bounds for arbitrary graphs was known in the bipartite setting~\cite{DHJKV}.

Let $\mathcal M(G)$ be the set of induced matchings of $G$. The
\emph{fractional strong chromatic index} of $G$ is
\[
\csf(G)=\min\Big\{\sum_{M\in\mathcal M(G)} w_M \;:\;
w_M\ge 0 \text{ for all } M,\
\sum_{M\in\mathcal M(G),\, M\ni e} w_M\ge 1 \text{ for every } e\in E(G)\Big\}.
\]
Requiring $w_M\in\{0,1\}$ gives $\cs(G)$, so $\csf(G)$ is the optimum of the
linear programming relaxation of this integer program, and
$\csf(G)\le\cs(G)$. Equivalently, $\csf(G)$ is the fractional chromatic
number of $L(G)^2$, whose independent sets are exactly the induced matchings
of $G$; see \'Sleszy\'nska-Nowak~\cite{SN}.
 Cames van Batenburg, Kang and
Pirot~\cite{CKP} showed that $\csf(G)\le1.5\,\Delta(G)^2$ for bipartite $G$;
D\k{e}bski~\cite{D} improved it to $1.4762\Delta^2$.

Our first result is the fractional version of Conjecture~\ref{conj:BQM}.

\begin{theorem}\label{thm:frac}
Every bipartite graph $G$ with parts $A$ and $B$ satisfies $\csf(G)\le\Delta_A\Delta_B$.
\end{theorem}

Since an induced matching in the complete bipartite graph $K_{\Delta_A, \Delta_B}$ contains exactly one edge, $\csf(K_{\Delta_A, \Delta_B} ) = \Delta_A\Delta_B$ and Theorem~\ref{thm:frac} is tight. Our second result is an asymptotic version of Conjecture~\ref{conj:BQM} with an explicit error term.

\begin{theorem}\label{thm:main}
Let $G$ be a bipartite graph with parts $A$ and $B$ and $\Delta=\max\{\Delta_A,\Delta_B\}$. Then
\[
\cs(G)\le(1+o(1))\,\Delta_A\Delta_B,
\]
where $o(1)$ is a quantity  depending only on  $\Delta$ that tends to $0$ as $\Delta\to\infty$. More precisely, 
\[
\cs(G)\le\Big(1+12\,\frac{(\ln\Delta)^{1/2}}{\Delta^{1/4}}\Big)\Delta_A\Delta_B .
\]
\end{theorem}





The main ingredient in our proofs is a simple coloring lemma based on linear orders of one part of the bipartition. Given a linear order of $A$, let each vertex $b\in B$ select the edge joining it to its first neighbor in that order.
Lemma~\ref{lem:column} shows that the selected edges can be strongly colored
with at most $\Delta_A$ colors. Choosing the order at random then yields the
fractional bound in Theorem~\ref{thm:frac}.

For the ordinary strong chromatic index, we use several linear orders of $A$
and their reverses. These orders produce a collection of columns to which the
edges are assigned. Lemma~\ref{lem:column} reduces the problem to balancing
the number of edges assigned to each column around each vertex of $A$. We
obtain the required balance using the Lov\'asz Local Lemma, which leads to
Theorem~\ref{thm:main}.

The paper is organized as follows. Section~\ref{sec:prelim} collects the probabilistic tools and proves the column lemma; both are key ingredients in our proof. The proof of Theorem~\ref{thm:frac} is presented in Section~\ref{sec:frac}. Finally, in Section~\ref{sec:main}, we first construct balanced families of orders and develop the randomized rounding argument, and then use them to prove Theorem~\ref{thm:main}.

\section{Preliminaries}\label{sec:prelim}

For a vertex $v$ we write $N(v)$ for its neighborhood and $d(v)=|N(v)|$ for its degree; $\ln$ is the natural logarithm and $[n]=\{1,\dots,n\}$. Throughout, $G$ is a bipartite graph with parts $A$ and $B$; we write each edge as $xb$ with $x\in A$ and $b\in B$. Vertices of degree $0$ play no role, and we assume that every vertex of $B$ has degree at least $1$.

Two distinct edges of a graph \emph{conflict} if they share an endpoint or are joined by an edge. Since $G$ is bipartite, two distinct edges $xb$ and $x'b'$, where $x,x'\in A$ and $b,b'\in B$, conflict if and only if $x=x'$, $b=b'$, $xb'\in E(G)$, or $x'b\in E(G)$. Hence a set of edges is an induced matching if and only if no two of its edges conflict, and a strong edge-coloring is a coloring of $E(G)$ in which conflicting edges receive distinct colors.

We use two probabilistic tools. The first is a form of the Chernoff and Bernstein inequalities, see \cite[Theorems~2.1 and~2.8]{JLR}.

\begin{lemma}\label{lem:chernoff}
Let $S$ be a sum of finitely many independent random variables with values in $\{0,1\}$, and let $\mu$ be the expected value of $S$. Then for every $t\ge0$,
\[
\Pr(S\ge\mu+t)\le\exp\Big(-\frac{t^2}{2(\mu+t/3)}\Big)\quad\text{and}\quad\Pr(S\le\mu-t)\le\exp\Big(-\frac{t^2}{2\mu}\Big).
\]
\end{lemma}

The second is the symmetric Lov\'asz Local Lemma, see \cite[Corollary~5.1.2]{AS}.

\begin{lemma}[Lovász Local Lemma]\label{lem:LLL}
Let $\cO_1,\dots,\cO_n$ be events, each having probability at most $p$.
Suppose that each $\cO_i$ is mutually independent of all but at most $d$
of the other events. If $e \cdot p \cdot (d+1) \le 1$,  then $\Pr\left(\bigcap_{i=1}^n \overline{\cO_i}\right)>0$.
In particular, with positive probability none of the events
$\cO_1,\dots,\cO_n$ occurs.
\end{lemma}


For our purposes, the underlying probability space is a product of finitely many independent random variables, and each event $\cO_i$ is determined by a set $S_i$ of these variables. In this situation $\cO_i$ is mutually independent of the family of all events $\cO_j$ with $S_j\cap S_i=\emptyset$, which is how the hypothesis of Lemma~\ref{lem:LLL} will be verified.


 Let $K >0$ be an integer and let $\sigma_1,\dots,\sigma_K$ be linear orders of $A$. For each $i \in [K]$, associate two columns, denoted by $(i,-)$ and $(i,+)$, corresponding to rules to assign the first and last neighbors of $b\in B$ in the order $\sigma_i$, respectively. Let $J = \{(i,-),(i,+): i \in[K]\}$, and for each $b \in B$ and $j \in J$, the {\em center} $c_j(b) \in N(b)$ is defined by
\[
c_j(b)=\begin{cases}
\min_{\sigma_i}N(b),\qquad  j=(i,-)\\
\max_{\sigma_i}N(b),\qquad j=(i,+).
\end{cases}
\]
For $y\in N(b)$, let $M_y(b)=\{j\in J:c_j(b)=y\}$ be the set of columns in which $y$ is the center of $b$, and let $m_y(b)=|M_y(b)|$. For fixed $b$, $M_y(b) \cap M_{y'}(b) = \emptyset$, where $y \neq y' \in N(b)$,  and $\cup _{y \in N(b)} M_y(b) = J$. If $d(b)\ge2$, then $m_y(b)$ is the number of indices $i$ for which $y$ is the $\sigma_i$-smallest or the $\sigma_i$-largest vertex of $N(b)$; if $d(b)=1$, then $m_y(b)=2K$.

\begin{lemma}[Column lemma]\label{lem:column}
Let $G$ be a bipartite graph with parts $A$ and $B$. 
Let $V_j\subseteq B$ for $j\in J$, and let $E_j=\{c_j(b)b:b\in V_j\}$.
\begin{enumerate}
\item[(a)] For every $j\in J$, the edge set $E_j$ can be partitioned into at most $\max_{x\in A}|N(x)\cap V_j|$ induced matchings of $G$.
\item[(b)] If every edge of $G$ lies in at least one of the sets $E_j$, then
\[
\cs(G)\le\sum_{j\in J}\max_{x\in A}|N(x)\cap V_j| .
\]
\end{enumerate}
\end{lemma}

\begin{proof}
(a) Let $j=(i,-)$; for $j=(i,+)$ replace $\sigma_i$ by its reverse in the following argument. Put $\lambda=\max_{x\in A}|N(x)\cap V_j|$. Every $b\in V_j$ contributes exactly one edge to $E_j$, namely $c_j(b)b$. Color the edges of $E_j$ one by one, processing their centers in increasing $\sigma_i$-order and the edges with a common center in any order, and give each edge a color in $[\lambda]$ that is not used on an already colored edge conflicting with it.

We show that such a color always exists. Let $xb\in E_j$, so $x=c_j(b)$, and let $x'b'\in E_j$ be an already colored edge that conflicts with $xb$. Then $b'\ne b$, since $b$ contributes only one edge to $E_j$, and $x'\le_{\sigma_i}x$, since centers are processed in increasing order. If $x'\ne x$ and $x'b\in E(G)$, then $x'$ is a neighbor of $b$ with $x'<_{\sigma_i}x=\min_{\sigma_i}N(b)$, which is impossible. As $b\ne b'$, the definition of conflict leaves two possibilities: $x'=x$, or $xb'\in E(G)$. In both cases $x b' \in E(G)$, and hence, $b'\in N(x)$. Therefore, every already colored edge conflicting with $xb$ is of the form $c_j(b')b'$ with $b'\in(N(x)\cap V_j)\setminus\{b\}$. There are at most $\lambda-1$ such edges, so a color in $[\lambda]$ is available. In the resulting coloring no two conflicting edges have the same color, so every color class is an induced matching of $G$.

(b) Assign each edge to one of the sets $E_j$ containing it. Color the edges assigned to $E_j$ as in (a), using a palette of $\max_x|N(x)\cap V_j|$ colors, with pairwise disjoint palettes for different columns. 
Since distinct columns use disjoint palettes, any two edges receiving the same color belong to the same $E_j$, and hence do not conflict by part (a).
\end{proof}

\begin{figure}[t]
\centering
\begin{tikzpicture}[xscale=1.25,yscale=1.1,every node/.style={font=\small}]
\foreach \i/\x in {1/0,2/1.2,3/2.4,4/3.6,5/4.8}{\node[circle,fill,inner sep=1.5pt,label=below:$x_{\i}$] (x\i) at (\x,0) {};}
\foreach \i/\x in {1/0.6,2/1.8,3/3.0,4/4.2}{\node[circle,draw,fill=white,inner sep=1.5pt,label=above:$b_{\i}$] (b\i) at (\x,1.8) {};}
\foreach \a/\b in {x3/b1,x3/b2,x5/b2,x4/b3,x4/b4,x5/b4}{\draw[gray] (\a)--(\b);}
\draw[very thick] (x1)--(b1); \draw[very thick] (x2)--(b2);
\draw[very thick,dashed] (x1)--(b4); \draw[very thick,dashed] (x3)--(b3);
\node at (6.2,0) {$A$}; \node at (6.2,1.8) {$B$};
\end{tikzpicture}
\caption{The order $\sigma$ of $A$ is $x_1<x_2<x_3<x_4<x_5$. Each $b_i$ keeps the edge to its $\sigma$-first neighbor (thick edges). Lemma~\ref{lem:column}(a), with one column and $V_j=B$, guarantees a strong coloring of the thick edges with $\max_x d(x)=3$ colors; here the procedure of its proof, always taking the smallest available color, uses two: solid edges $x_1b_1,x_2b_2$ and dashed edges $x_1b_4,x_3b_3$. Note that $x_3b_3$ and $x_1b_1$ conflict (via the edge $x_3b_1$), as do $x_3b_3$ and $x_2b_2$ (via $x_3b_2$), while $x_3b_3$ and $x_1b_4$ do not.}
\label{fig:column}
\end{figure}
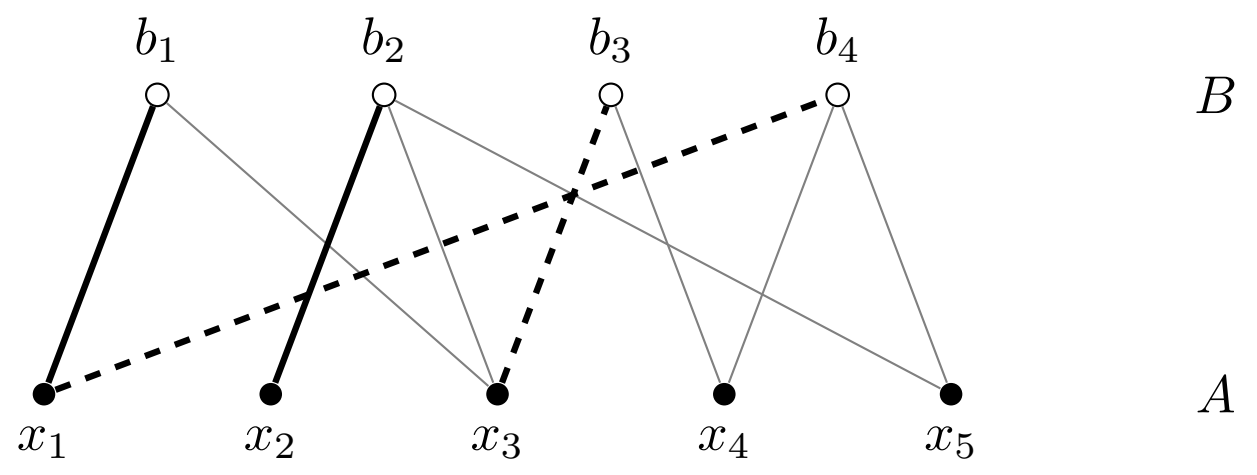

\begin{remark}
The choice of the first neighbor is used exactly once in the proof: it excludes conflicts through an edge $x'b$ with $x'$ preceding $x$. Figure~\ref{fig:column} shows a small example. The bound in Lemma~\ref{lem:column}(a) is attained, for instance, by a complete bipartite graph with a single column and $V_j=B$.
\end{remark}

\section{proof of Theorem~\ref{thm:frac}: the fractional bound}\label{sec:frac}

The proof of Theorem~\ref{thm:frac} has two parts. First, a fractional strong edge-coloring can be obtained from a suitable
probability distribution on the induced matchings of $G$
(Lemma~\ref{lem:randmatch}). Second, a single random order of $A$ produces a random induced matching that contains every edge with probability at least $1/(\Delta_A\Delta_B)$.

\begin{lemma}\label{lem:randmatch}
Suppose there exists a probability distribution on the induced matchings of $G$ such that a random induced matching $M$ drawn from this distribution satisfies $\Pr(e\in M)\ge p$ for every edge $e\in E(G)$. Then $\chi'_{s,f}(G)\le \frac{1}{p}$.
\end{lemma}

\begin{proof}
For every induced matching $M'$ of $G$ put $w_{M'}=\Pr(M=M')/p$. These weights are nonnegative, and since the probabilities $\Pr(M=M')$ add up to $1$, the weights add up to $1/p$. For every edge $e$,
\[
\sum_{e \in M'}w_{M'}=\frac1p\sum_{e \in M'}\Pr(M=M')=\frac{\Pr(e\in M)}{p}\ge1 .
\]
So the weights $w_{M'}$ are admissible in the definition of $\csf(G)$, and $\csf(G)\le\sum_{M'}w_{M'}=1/p$.
\end{proof}

\begin{proof}[Proof of Theorem~\ref{thm:frac}]
By Lemma~\ref{lem:randmatch}, it suffices to construct a probability
distribution on the induced matchings of $G$ such that a random induced
matching $M$ drawn from this distribution satisfies
\[
\Pr(e\in M)\ge \frac{1}{\Delta_A\Delta_B}
\]
for every edge $e\in E(G)$. We construct such a distribution in four steps.

\smallskip\noindent\emph{Step 1: a random order.} Let $\sigma$ be a linear order of $A$ chosen uniformly at random from all linear orders of $A$.

\smallskip\noindent\emph{Step 2: every vertex of $B$ keeps one edge.} Every $b\in B$ keeps the edge that joins $b$ to its $\sigma$-first neighbor. Let
\[
E_\sigma=\big\{(\min\nolimits_\sigma N(b))\,b:\ b\in B\big\}
\]
be the set of kept edges; it contains exactly one edge at each vertex of $B$.

\smallskip\noindent\emph{Step 3: splitting the kept edges.} Apply Lemma~\ref{lem:column}(a) with $K=1$, $\sigma_1=\sigma$, the column $j=(1,-)$ and $V_j=B$. Then $E_j=E_\sigma$ and $\max_{x\in A}|N(x)\cap V_j|=\max_{x\in A}d(x)=\Delta_A$. Hence $E_\sigma$ can be partitioned into induced matchings $M_{\sigma,1},\dots,M_{\sigma,\Delta_A}$ of $G$, some of which may be empty. We fix such a partition for every linear order $\sigma$ of $A$.

\smallskip\noindent\emph{Step 4: choosing one part.} Choose $t\in[\Delta_A]$ uniformly at random and independently of $\sigma$, and put $M=M_{\sigma,t}$. Then $M$ is an induced matching of $G$.

\smallskip
Now fix an edge $yb$ with $y\in A$ and $b\in B$. The edge $yb$ belongs to $M$ exactly when the following two events occur:
\begin{enumerate}
\item[(i)] $yb\in E_\sigma$, that is, $y$ is the $\sigma$-first vertex of $N(b)$;
\item[(ii)] $t$ is the index of the part of the partition of $E_\sigma$ that contains $yb$.
\end{enumerate}
The order that $\sigma$ induces on $N(b)$ is uniformly random, so each of the $d(b)$ vertices of $N(b)$ is the $\sigma$-first one with the same probability, and (i) has probability $1/d(b)$. If (i) occurs, then $yb$ lies in exactly one of the parts $M_{\sigma,1},\dots,M_{\sigma,\Delta_A}$, and since $t$ is uniform on $[\Delta_A]$ and independent of $\sigma$, the conditional probability of (ii) is $1/\Delta_A$. Therefore
\[
\Pr(yb\in M)=\frac{1}{d(b)}\cdot\frac{1}{\Delta_A}\ge\frac{1}{\Delta_A\Delta_B},
\]
since $d(b)\le\Delta_B$. Lemma~\ref{lem:randmatch} with $p=1/(\Delta_A\Delta_B)$ gives $\csf(G)\le\Delta_A\Delta_B$.
\end{proof}


\section{The proof of Theorem~\ref{thm:main}}\label{sec:main}

First, suppose that $k\ge2$ and $\Delta\ge1$ are integers, and $G$ is a bipartite graph with parts $A$ and $B$ such that $d(x)\le\Delta$ for every $x\in A$ and $1\le d(b)\le k$ for every $b\in B$. A vertex $y\in N(b)$ is the center of $b$ in exactly $m_y(b)$ of the $2K$ columns. On average this number is $2K/d(b)\ge2K/k$; the next lemma shows that all of them can be made close to this value simultaneously.

\begin{lemma}\label{lem:bal}
Let $\delta\in(0,1)$ and let $K$ be an integer with \[K\ge\frac{k}{\delta^2}\ln(3k^2\Delta).\] There are linear orders $\sigma_1,\dots,\sigma_K$ of $A$ such that
\[
m_y(b)\ge(1-\delta)\,\frac{2K}{k}\qquad\text{for every }b\in B\text{ and every }y\in N(b).
\]
\end{lemma}

\begin{proof}
If $d(b)=1$, then $m_y(b)=2K$. If $d(b)=2$, then in every order one vertex of $N(b)$ is the smallest and the other is the largest, so $m_y(b)=K\ge2K/k$. These vertices $b$ satisfy the conclusion for any orders.

Choose independent random labels $\ell_i(x)$, uniform in $[0,1]$, for $i\in[K]$ and $x\in A$, and let $\sigma_i$ order $A$ by increasing $\ell_i$; ties occur with probability $0$. Let $b\in B$ with $d=d(b)\ge3$ and let $y\in N(b)$. For each $i$, the events ``$y$ is the $\sigma_i$-smallest vertex of $N(b)$'' and ``$y$ is the $\sigma_i$-largest vertex of $N(b)$'' are disjoint  with  probability $1/d$ each, and those events with different $i$s are independent. So $m_y(b)$ is a sum of $K$ independent indicator variables with mean $\mu=2K/d\ge2K/k$. By Lemma~\ref{lem:chernoff} with $t=\delta\mu$,
\begin{align*}
\Pr\Big(m_y(b)<(1-\delta)\frac{2K}{k}\Big)&\le\Pr\big(m_y(b)\le(1-\delta)\mu\big)
\le\exp\Big(-\frac{\delta^2\mu}{2}\Big)\le\exp\Big(-\frac{\delta^2K}{k}\Big).
\end{align*}
Let $\cD_b$ be the event that $m_y(b)<(1-\delta)2K/k$ for some $y\in N(b)$. Then $\Pr(\cD_b)\le k\exp(-\delta^2K/k)=:p$. The event $\cD_b$ is determined by the labels of the vertices in $N(b)$, so it is mutually independent of the family of events $\cD_{b'}$ with $N(b')\cap N(b)=\emptyset$. The number of vertices $b'\ne b$ with $N(b')\cap N(b)\ne\emptyset$ is at most $\sum_{x\in N(b)}(d(x)-1)\le k(\Delta-1)$. Since $K\ge(k/\delta^2)\ln(3k^2\Delta)$,
\[
e\cdot p\cdot\big(k(\Delta-1)+1\big)\le e\cdot k\cdot k\Delta\cdot\exp\Big(-\frac{\delta^2K}{k}\Big)\le\frac{e}{3}<1 .
\]
By Lemma~\ref{lem:LLL}, with positive probability none of the events $\cD_b$ occurs.
\end{proof}

The next result may be regarded as a version of Theorem~\ref{thm:main} with free parameters.

\begin{theorem}\label{thm:rand}
Let $\delta\in(0,1)$, let $K$ be an integer with $K\ge(k/\delta^2)\ln(3k^2\Delta)$, and put $L=\ln(9Kk\Delta)$. Then
\[
\cs(G)\le\frac{k\Delta}{1-\delta}+2\sqrt{\frac{Kk\Delta L}{1-\delta}}+2KL .
\]
\end{theorem}

\begin{proof}
Fix orders $\sigma_1,\dots,\sigma_K$ as in Lemma~\ref{lem:bal}. For every $b\in B$ and every $y\in N(b)$, $m_y(b) \geq (1- \delta) \frac{2K} {k} > 0$ by Lemma~\ref{lem:bal}. Hence $M_y(b) \neq \emptyset$. Choose a column $j_{b,y}\in M_y(b)$ uniformly at random and all choices are independent. For $j\in J$ put
\[
V_j=\{b\in B: j_{b,c_j(b)}=j\},
\]
and let $E_j=\{c_j(b)b:b\in V_j\}$ as in Lemma~\ref{lem:column}. Since $j_{b,y} \in M_y(b)$, every edge $yb$ lies in $E_{j}$ with $j=j_{b,y}$. 

Fix $x\in A$ and $j\in J$, and let
\[
X_{x,j}=|N(x)\cap V_j|=\sum_{b\in N(x)}\mathbf 1\big[j_{b,c_j(b)}=j\big].
\]
The summand corresponding to $b$ depends only on the choice
$j_{b,c_j(b)}$, and the choices corresponding to different vertices
$b\in N(x)$ are made independently. Since $j\in M_{c_j(b)}(b)$, the summand of $b$ equals $1$ with probability $1/m_{c_j(b)}(b)$, which is at most $k/(2(1-\delta)K)$ by Lemma~\ref{lem:bal} (and equals $1/(2K)$ if $d(b)=1$). As $d(x)\le\Delta$, the expected value $\mu_{x,j}$ of $X_{x,j}$ is at most
\[
\lambda:=\frac{k\Delta}{2(1-\delta)K}.
\]
If $\mu_{x,j}<\lambda$, add to $X_{x,j}$ the sum of $r=\lceil\lambda-\mu_{x,j}\rceil$ indicator variables, each with mean $(\lambda-\mu_{x,j})/r$, independent of each other and of all choices $j_{b,y}$. The result $X'_{x,j}$ is a sum of independent indicators with $X'_{x,j}\ge X_{x,j}$ and expected value exactly $\lambda$; if $\mu_{x,j}=\lambda$, put $X'_{x,j}=X_{x,j}$. Let $t=\sqrt{2\lambda L}+L$. Then
\[
t^2-2\lambda L-\tfrac23Lt=\tfrac43L\sqrt{2\lambda L}+\tfrac13L^2\ge0,
\]
that is, $t^2/(2(\lambda+t/3))\ge L$. By Lemma~\ref{lem:chernoff},
\[
\Pr(X_{x,j}\ge\lambda+t)\le\Pr(X'_{x,j}\ge\lambda+t)\le\exp(-L)=\frac{1}{9Kk\Delta}.
\]

Let $\cO_{x,j}$ be the event $X_{x,j}\ge\lambda+t$. It is determined by the choices $j_{b,y}$ with $b\in N(x)$, so it is mutually independent of the family of events $\cO_{x',j'}$ with $N(x')\cap N(x)=\emptyset$. The number of vertices $x'$ with $N(x')\cap N(x)\ne\emptyset$, including $x$ itself, is at most $1+\sum_{b\in N(x)}(d(b)-1)\le1+\Delta(k-1)\le k\Delta$; hence the number of pairs $(x',j')$ with $N(x')\cap N(x)\ne\emptyset$, including $(x,j)$, is at most $2Kk\Delta$. Since
\[
e\cdot\frac{1}{9Kk\Delta}\cdot2Kk\Delta=\frac{2e}{9}<1,
\]
Lemma~\ref{lem:LLL} shows that for some choice of the columns $j_{b,y}$ we have $X_{x,j}<\lambda+t$ for all $x\in A$ and $j\in J$. For this choice, Lemma~\ref{lem:column}(b) gives
\[
\cs(G)\le\sum_{j\in J}\max_{x\in A}X_{x,j}\le2K(\lambda+t)=\frac{k\Delta}{1-\delta}+2K\sqrt{2\lambda L}+2KL,
\]
and $2K\sqrt{2\lambda L}=2\sqrt{Kk\Delta L/(1-\delta)}$.
\end{proof}

Now we are ready to prove Theorem~\ref{thm:main}

\begin{proof}[Proof of Theorem~\ref{thm:main}]
Let
\[
u=\frac{(\ln\Delta)^{1/2}}{\Delta^{1/4}},\qquad\text{so that}\qquad u^2=\frac{\ln\Delta}{\Delta^{1/2}} .
\]
Choose $\Delta$ large enough so that $u<1/10$. 
Put
\[
\delta=u,\qquad K=\Big\lceil\frac{k}{u^2}\ln(3k^2\Delta)\Big\rceil,\qquad L=\ln(9Kk\Delta).
\]
Then $\delta\in(0,1)$ and $K$ satisfies the hypothesis of Theorem~\ref{thm:rand}. Dividing the bound of Theorem~\ref{thm:rand} by $k\Delta=\Delta_A\Delta_B$ gives
\begin{equation}\label{eq:threeterms}
\frac{\cs(G)}{\Delta_A\Delta_B}\le\frac1{1-u}+2\sqrt{\frac{1}{1-u}\cdot\frac{KL}{k\Delta}}+2\,\frac{KL}{k\Delta}.
\end{equation}

We estimate $K$ and $L$. Since $k\le\Delta$, we have $\ln(3k^2\Delta)\le\ln3+3\ln\Delta\le4\ln\Delta$. So 
\[
K=\Big\lceil\frac{k}{u^2}\ln(3k^2\Delta)\Big\rceil\le \Big\lceil\frac{4k\ln\Delta}{u^2}\Big\rceil=\Big\lceil4k\Delta^{1/2}\Big\rceil\le 4k\Delta^{1/2}+1.
\]
Hence $L=\ln (9Kk\Delta)\le \ln (9\cdot 4.001\cdot \Delta^{7/2})\le \ln 37+7/2\ln \Delta$. Therefore, 
\[
\frac{KL}{k\Delta}\le\frac{(4k\Delta^{1/2}+1)\cdot (\ln 37+7/2\ln \Delta)}{k\Delta}\le 4\cdot 3.6\,\frac{\ln\Delta}{\Delta^{1/2}}=14.4u^2\le 1.44\,u .
\]

Finally, $u\le1/10$ gives $\frac1{1-u}=1+\frac{u}{1-u}\le1+\frac{10}{9}u$ and $\frac1{1-u}\le\frac{10}{9}$. So the three terms on the right-hand side of \eqref{eq:threeterms} satisfy
\[
\frac1{1-u}\le1+\frac{10}{9}u,\qquad
2\sqrt{\frac{1}{1-u}\cdot\frac{KL}{k\Delta}}\le2\sqrt{\frac{10}{9}\cdot 14.4u^2}=8\,u,\qquad
2\,\frac{KL}{k\Delta}\le 2.88\,u .
\]
Adding them,
\[
\frac{\cs(G)}{\Delta_A\Delta_B}\le1+\Big(\frac{10}{9}+8+2.88\Big)u<1+12u .\qedhere
\]
\end{proof}

\begin{remark}
The choice of the side $A$ in the proof matters: the parameters $L_0$ and $L$ grow with $\ln k$, and the error term is small because $k\le\Delta$. The proof also shows where the error term comes from. The factor $1/(1-\delta)$ is the price for the imbalance of the counts $m_y(b)$ among $K$ random orders, which forces $K$ to be of order $k\log\Delta/\delta^2$; the term $2\sqrt{Kk\Delta L/(1-\delta)}$ is the price for the fluctuations of the loads $X_{x,j}$, whose expected values $\lambda$ are of order $k\Delta/K$. Balancing the two gives the relative error of order $(\log\Delta)^{1/2}\Delta^{-1/4}$. We could improve $12$ in Theorem~\ref{thm:main}  to a smaller constant if we are more careful in calculation. 
\end{remark}

\section*{AI usage declaration}
Claude Opus~5.5 was used to assist in developing the proofs and in preparing the text of this article. The main ideas are due to the authors, who have verified all mathematical content and take full responsibility for its correctness.

\end{document}